\documentclass[24pt]{article}
\usepackage{amssymb,amsmath,makeidx,amscd,verbatim,amsthm,graphicx,mathrsfs}
\usepackage{lscape}
\usepackage{color}
\newtheorem{thm}{Theorem}[section]

\newtheorem{lem}[thm]{Lemma}
\newtheorem{coro}[thm]{Corollary}
\theoremstyle{definition}
\newtheorem{defx}[thm]{Definition}

\numberwithin{equation}{section}

\newcommand{\be}{\begin{enumerate}}
	\newcommand{\ee}{\end{enumerate}}
\newcommand{\bq}{\begin{eqnarray*}}
	\newcommand{\eq}{\end{eqnarray*}}

\begin{document}
	\pagenumbering{arabic} \baselineskip 24pt
	\newcommand{\disp}{\displaystyle}
	\renewcommand*\contentsname{Table of Contents}
	\thispagestyle{empty}
	%\chapter{}
	%\chapter{}
	%\chapter{}
	\newcommand{\HRule}{\rule{\linewidth}{0.1mm}}
	\linespread{1.0}
	\pagenumbering{arabic} \baselineskip 10pt
	\thispagestyle{empty}
	\title{\textbf{Grushin Operator on Infinite Dimensional Homogeneous Lie Grougs}}
	\author{M. E. Egwe$^1$ and J. I. Opadara$^2$\\ Department of Mathematics, University of Ibadan, Ibadan.\\  $^1$\emph{murphy.egwe@ui.edu.ng},\;$^2$\emph{jossylib@gmail.com}}
	\maketitle
	\large	
	\begin{abstract}
		\noindent
		\textit{A collection of infinite dimensional complete vector fields $\left\{V_i\right\}_{i=1}^{\infty}$ acting on a locally convex manifolds $M$ on which a smooth positive measure $\mu$ is defined was considered. It was assumed that the vector fields generates an infinite dimensional Lie algebra $\mathfrak{g}$ and satisfies H$\ddot{o}$rmander's condition. The sum of squares of Grushin operators related to the vector fields was  examined and the operator is then considered as the generalized Grushin operator. The paramount proofs were Poincar$\acute{e}$ inequality, Gaussian two-bounded estimate for the related heat kernels and the doubling condition for the metric defined by the underlying vector fields.}
	\end{abstract}
	Keyword:Grushin Operator, doubling condition, locally convex manifold, H$\ddot{o}$rmander's condition.\\
	MSC Classification: 22Exx, 47AXX, 47DXX, 32W05
	
	\pagenumbering{arabic} \baselineskip 24pt
	\section{Introduction}
	Throughout this paper,  $G$ represents an infinite dimensional homogeneous Lie group and $\mathfrak{g}$ the Lie algebra of $G.$ The set $\mathbb{K}$ shall denote the set  of real numbers, $\mathbb{R}$ or set of complex numbers, $\mathbb{C}$.\\
	The Grushin operator is defined by the formula
	\begin{equation}\label{aaA}
	\mathcal{L}_n=-(\partial_1^2-v_1^{2n}\partial_2^2)
	\end{equation}
	on the space of $L^2(\mathbb{R}^2)$ where $n\in\mathbb{N}$ . Many scholars have used the operators in study of Lie groups and generalized it in many ways (cf. $\cite{2},\; \cite{4}, \; \cite{13}, \; \cite{15}\; \mbox{and}\; \cite{16}$).\\
	Equation$\eqref{aaA}$ above can be written as a sum of squares of vector fields
	\begin{equation*}
	\mathcal{L}=-V_1^2-V_2^2
	\end{equation*}
	where $V_1=\partial_1$ and $V_2=v_2^n\partial_2$ are two vector fields on $\mathbb{R}^2$. Recall that the vector fields $\left\{V_i\right\}_{i=1}^2$ span a finite dimensional Lie algebra $\mathfrak{g}$. Of course, this gives a starting point of this work.\\
	The main contribution of this paper is to show that in the context of simply connected Lie groups, doubling condition and Poincar$\acute{e}$ inequality, defined in $\eqref{KKK}$ and $\eqref{NNN}$ holds automatically.

	\subsection{Motivation and Setting of the Work}
	Notable motivations for this work were from the results of Jerison ($\cite{11}$); Jerison and S$\acute{a}$nchez-Calle ($\cite{12}$) and Dziuba$\acute{n}$ski and Sikora ($\cite{3}$). Jerison and S$\acute{a}$nchez-Calle proved that Poincar$\acute{e}$ inequality for vector fields satisfy H$\ddot{o}$rmander's condition and local bounds for heat kernels relating to the sum of squares of the underlined vector fields, while Dziuba$\acute{n}$ski and Sikora proved that the global Poincar$\acute{e}$ inequality and global heat kernels still hold in the setting of finite dimensional Lie groups. Here, the results can simply be stated as follows: if in addition that the vector fields are of infinite dimension which generate infinite dimensional Lie algebra $\mathfrak{g}$, the global Poincar$\acute{e}$ and heat kernel bounds also hold and satisfy the H$\ddot{o}$rmander's condition.\\
	Now, recall that a simply connected Lie group $G$ has  a polynomial growth if for all $V\in \mathcal{L}(G)=\mathfrak{g}$, the operator ad($V$) (adjoint of $V$) has only purely imaginary eigenvalues. Here,  doubling condition, the Gaussian two-sided bounds corresponding to heat kernels shall be proved for Grushin operators such as
	\begin{equation*}
	-(\partial_1^2+\partial_2^2+(v_1^2+v_2^2-1)^2\partial_3^2+\cdots)=-(V_1^2+V_2^2+V_3^2+\cdots)
	\end{equation*}
	where $V_1=\partial_1,\;\; V_2=\partial_2,\;\; V_3=(v_1^2+v_2^2-1)\partial_3$ and so on.\\
	Another good illustration is
	\begin{equation*}
	-(\partial_1^2+\sin^2v_1\partial_2^2+\cos^2v_1\partial_3^2+\cdots)=-(V_1^2+V_2^2+V_3^2+\cdots)
	\end{equation*}
	where $V_1=\partial_1,\;\; V_2=\sin v_1\partial_2,\;\; V_3=\cos v_1\partial_3$ and so on.

	\section{Preliminaries}
	Here we shall give the definition of an infinite-dimensional Lie group and give explanation of how its Lie algebra could be defined so as to define a functor from the class of Lie groups to the class of locally convex Lie algebras.
	\begin{defx}
		\ \\
		\normalfont
		(a) Let $V$ be a $\mathbb{K}$-vector space. A seminorm is a function $p:V\longrightarrow \mathbb{R}^+$ which satisfies the following axioms:\\
		(i) $p(\lambda v)=|\lambda|p(v)\;\;\;\; \forall\;\;\; \lambda\in \mathbb{K}, \; v\in V$\\
		(ii) $p(v+u)\leq p(v)+p(u)\;\;\;\; \forall\;\; u,v\in V$.\\
		In addition, if $p(v)>0\;\;\;\; \forall\;\; v\in V\setminus\{0\}$, then $p$ is called a norm.\\
		As a notation, the family of seminorms shall be denoted by $\mathcal{P}$.\\
		\\
		(b) The set $\mathcal{P}$ of seminorms $p$ on $V$ is called separating if $p(u)=0\;\;\; \forall\;\; p\in \mathcal{P}$ implies $u=0$.\\
		This is equivalent to the bijective linear map
		\begin{equation*}
		\beta:V\longrightarrow\underset{p\in\mathcal{P}}{\prod}V_p:u\mapsto(\beta(u))_{p\in\mathcal{P}}
		\end{equation*}
	\end{defx}
	\begin{defx}
		\ \\
		\normalfont
		A vector space $V$ is said to be locally convex  if it is endowed with a topology defined by a separating family of seminorms.
	\end{defx}
	\noindent
	Note that each locally convex space is in particular a topological vector space which can be embedded into a product $\underset{p\in\mathcal{P}}{\prod}V_p$ of normed spaces.
	\begin{defx}
		\ \\
		\normalfont
		If $V$ is a vector space and the map $\beta_i:V_i\longrightarrow V$ linear maps, defined on locally convex spaces $V_i$. We consider the system $\mathcal{P}$ of all the seminorms $p$ on $V$ for which all  compositions $p\circ\beta_i$ are continuous seminorms on the spaces $V_i$. By $\mathcal{P}$, one can obtain on $V$ a locally convex topology called the finally locally convex topology which is defined by the mappings $(\beta_i)_{i\in J}$.
	\end{defx}
	\begin{defx}
		\ \\
		\normalfont
		Let $G$ be a Lie group. Then for each $g\in G$, the map
		\begin{equation*}
		c:G\longrightarrow G:v\mapsto gvg^{-1}
		\end{equation*}
		is a smooth automorphism, hence prompts a continuous linear automorphism.
		\begin{equation*}
		Ad(g):=\mathcal{L}(c): \mathfrak{g}\longrightarrow\mathfrak{g}
		\end{equation*}
		We thus obtain an action $G\times \mathfrak{g}\longrightarrow\mathfrak{g}:(g,v)\mapsto Ad(g).v$ called adjoint action of $G$ on $\mathfrak{g}$
	\end{defx}
	\noindent
	Before proceeding to the definition of infinite dimensional Lie groups, the following notations are very useful.\\
	If $G$ is a Lie group and $g\in G$, we say that for all $v,w \in G$
	\begin{equation*}
	\beta:G\longrightarrow G: v\mapsto gv\;\;\; \mbox{is the left multiplicatiom by}\; g
	\end{equation*}
	\begin{equation*}
	\beta:G\longrightarrow G: v\mapsto vg\;\;\; \mbox{is the right multiplicatiom by}\; g
	\end{equation*}
	\begin{equation}\label{Ss}
	\beta:G\times G\longrightarrow G: (v,w)\mapsto vw\;\;\; \mbox{is the multiplicatiom map and}
	\end{equation}
	\begin{equation}\label{Tt}
	\beta:G\longrightarrow G: v\mapsto v^{-1}\;\;\; \mbox{is the inverse map}
	\end{equation}
	\begin{defx}
		\ \\
		\normalfont
		A locally convex manifold which is endowed with a group structure such that $\eqref{Ss}$ and $\eqref{Tt}$ are smooth is called a locally convex Lie group.\\
		Hence, an infinite dimensional Lie group $G$ is a group with an infinite dimensional manifold such that the two structures  are compatible, which means that  the multiplication map $\eqref{Ss}$ and inversion map $\eqref{Tt}$ are smooth.
	\end{defx}
	\noindent
	Such a Lie group $G$ is locally diffeomorphic to an infinite dimensional vector space $V_i\;\;\;i=1,2,3\cdots$, which could be a Banach space $\mathcal{B}$ induced by the norm $\parallel \cdot\parallel$, a Hilbert space $\mathcal{H}$ induced by inner product $\langle\cdot,\cdot\rangle$ and norm $\parallel x\parallel^2=\langle x,x\rangle$, or a Frechet space $\mathcal{F}$ induced by metric $\rho(\cdot, \cdot)$ but not norm. Correspondingly, we call these Banach Lie groups, Hilbert Lie groups or Frechet Lie groups.

	\section{Main Result}
	Suppose that $M$ is a locally convex manifold of infinite dimension with a positive measure $\mu$. It is assumed that the measure $\mu$ has a locally convex density with respect to any coordinate map on $M$. The tangent bundle of $M$ shall be denoted by $TM$ whose sections are vector fields on $M$. Let us consider an infinite family of vector fields $\left\{V_i\right\}_{i=1}^{\infty}$. It will be assumed that the vector fields $\left\{V_i\right\}_{i=1}^{\infty}$ spanned the map $\exp(sV_i)$ and is defined globally on the whole of $M$, for all real number $s$ and $i=1,2,\cdots$. It is sufficient to say that $V_i$ are complete vector fields. We lay down some assumptions.\\
	(1). Recall that if $U$ and $V$ are any two vector fields, their commutator which is also a vector field is defined as
	\begin{equation*}
	[U,V]f=UVf-VUf.
	\end{equation*}
	Consequently, it will also be assumed that the collection $\{V_1,V_2, \cdots\}$ in conjunction with all their commutators span an infinite dimensional Lie algebra $\mathfrak{g}$.\\
	(2). In the sequel, it will be assumed that the vector fields $\left\{V_i\right\}_{i=1}^{\infty}$ are skew-adjoint, that is
	\begin{equation}\label{aa}
	\displaystyle\int_MV_if(v)g(v)d\mu(v)=-\displaystyle\int_MV_if(v)g(v)d\mu(v).
	\end{equation}
	In short one can write $V_i^{\ast}=-V_i$ for $\eqref{aa}$\\
	Now, if $U$ and $V$ are skew-adjoint and $W=[U,V]$, related to $\eqref{aa}$ above, we have that $W^{\ast}=-W$. It is worthy to note that, all the elements of a Lie algebra $\mathfrak{g}$ generated by skew-adjoint vector fields are skew-adjoint by implication..\\
	For a complete skew-adjoint vector field $V$ on manifold $M$ we have
	\begin{equation*}
	\dfrac{d}{ds}\displaystyle\int_M|f(\exp(sV)v)|^2d\mu(v)=\displaystyle\int_M[(V+V^{\ast})f(\exp(sU)v)]\overline{f(\exp(sV)v)d\mu(v)}=0,
	\end{equation*}
	where $f\in C_{loc}^{\infty}(M)$. Therefore, the graph $\exp(sV)$ preserves the positive measure $\mu$, which means
	\begin{equation}\label{JJJ}
	\displaystyle\int|f(v)|^pd\mu(v)=\displaystyle\int|f(\exp(sV)v)|^pd\mu(v)
	\end{equation}
	$\forall \; V\in \mathfrak{g},\;\; f\in L^p(M)$ being an integrable function, where $p\in[1,\infty]$ and $s$ is as defined before.\\
	If the infinite number of commutators of $\{V_1,V_2,\cdots\}$ linearly spans the tangent space $T_vM$, then we say that the collection of vector fields $\{V_1,V_2,\cdots\}$ satisfy H$\ddot{o}$rmander condition, for all $v \in M$.\\
	Since it was assumed that the collection $\left\{V_i\right\}_{i=1}^{\infty}$ generates the infinite dimensional Lie algebra $\mathfrak{g},$ then, in this regard, H$\ddot{o}$rmander's condition implies, for any $v\in M$ the linear space corresponding to $\mathfrak{g}$ at $v$ is $TM_v$.\\
	Now, let us define the control distance or the sub-Riemannian distance between any two points $v,w \in M$ by $\rho(v,w)$ and the open ball with respect to $\rho$ by $B(v,r)$ with centre at $v$ and radius $r$. Then we have $U(v,r)=\mu(B(v,r))$. Hence a metric measure space satisfies the doubling condition when
	\begin{equation}\label{KKK}
	U(v,2r)\leq KU(v,r),\;\;\;\;\;v\in M\;\; \mbox{and}\;\; r>0.
	\end{equation}
	The focus in this setting is the sum of squares of operator corresponding to the collection of vector fields $\left\{V_i\right\}_{i=1}^{\infty}$ which by definition is
	\begin{equation}\label{LLL}
	\mathcal{L}=\displaystyle\sum_{i=1}^{\infty}V_iV_i^{\ast}=-\displaystyle\sum_{i=1}^{\infty}V_i^2.
	\end{equation}
	The equivalent gradient is defined for a function $f$ in $C_{loc}^{\infty}(M)$ by
	\begin{equation*}
	\bigtriangledown f=(V_1f,V_2f,V_3f,\cdots),
	\end{equation*}
	and setting
	\begin{equation*}
	|\bigtriangledown f(v)|^2=|V_1f+V_2f+V_3f+\cdots|^2=\displaystyle\sum_{i=1}^{\infty}|V_if(v)|^2,
	\end{equation*}
	by implication, we have
	\begin{equation*}
	\parallel \bigtriangledown f\parallel_2=\parallel\mathcal{L}^{1/2}f\parallel_2.
	\end{equation*}
	With the above facts, the main result can now be stated as follows:
	\begin{thm}\label{aAa}
		\ \\
		\normalfont
		Let $M$ be an infinite dimensional locally convex and smooth manifold and suppose that the infinite number of commutators of the smooth vector fields $\left\{V_i\right\}_{i=1}^{\infty}$ linearly spans $T_vM$ and generates an infinite dimensional Lie algebra $\mathfrak{g}$. Furthermore, assume that the vector fields $\left\{V_i\right\}_{i=1}^{\infty}$ are skew-adjoint.
		Then the distance $\rho$ satisfies $\eqref{KKK}$ and the heat kernel $h_s$ equivalent to the square operator $\mathcal{L}$ satisfies the two-sided Gaussian bounds
		\begin{equation}\label{MMM}
		\dfrac{ae^{-a\frac{d(u,v)^2}{s}}}{[(\mu(B(u,\sqrt{s})))(\mu(B(v,\sqrt{s})))]^{1/2}}\leq h_s(u,v)\leq \dfrac{a^{\prime}e^{-a^{\prime}\frac{d(u,v)^2}{s}}}{[(\mu(B(u,\sqrt{s})))(\mu(B(v,\sqrt{s})))]^{1/2}}.
		\end{equation}
		Moreso, for any $p\in(1,\infty],$ Riesz transformation is bounded, that is
		\begin{equation*}
		\parallel\bigtriangledown f\parallel_p\leq K_p\parallel\mathcal{L}^{1/2}f\parallel_p.
		\end{equation*}
		Note that the two-sided Gaussian estimate $\eqref{MMM}$ implies Poincar$\acute{e}$ inequality.
	\end{thm}
	\noindent
	Before proceeding to the proof of the above Theorem, the following Corollaries are essential
	\begin{coro}
		\ \\
		\normalfont
		Maintaining the assumption of the Theorem 2.1, the locally convex manifold $M$ satisfies the Poincar$\acute{e}$ Inequality, this means that $\exists$ a positive constant $K$ such that
		\begin{equation}\label{NNN}
		\displaystyle\int_B|f-f_B|^2d\mu\leq Kr^2\displaystyle\int_B|\bigtriangledown f|^2d\mu,
		\end{equation}
		where $B$ is any ball and $r$ is its radius
	\end{coro}
	\noindent
	The operator $\mathcal{L}$ defined in $\eqref{LLL}$ is positive and self-adjoint. Hence, the operator admits a spectral resolution $D(\gamma)$. and for any Borel function $H:[0,\infty]\longrightarrow \mathbb{C}$, the operator $H(\mathcal{L})$ can be defined as
	\begin{equation}\label{OOO}
	H(\mathcal{L})=\displaystyle\int_0^{\infty} H(\gamma)dD_{\mathcal{L}}(\gamma).
	\end{equation}
	Now, on the metric measure space, the doubling condition defined in $\eqref{KKK}$ means that there are constants $K,\nu$ such that for any $\gamma\geq1$ and $v\in V$ one has
	\begin{equation}\label{PPP}
	U(v,\gamma r)\leq K\gamma^{\nu}U(v,r).
	\end{equation}
	Now, it is believed that Gaussian bounds implies the following spectral multiplier which is another implication of Theorem 3.1.
	\begin{coro}
		\ \\
		\normalfont
		Suppose that $\nu$ in relation $\eqref{PPP}$ is an exponent corresponding to the locally convex manifold $M$ and the control distance $\rho$.\\
		Let $\dfrac{\nu}{2}<l$ and $H:[0,\infty]\longrightarrow \mathbb{C}$ be a Borel function such that
		\begin{equation}\label{QQQ}
		\underset{0<s}{\sup}\parallel\alpha\delta_sH\parallel_{E^{l,\infty}}\leq\infty
		\end{equation}
		where $\delta_sH(\gamma)=H(s\gamma)$ and $\parallel H\parallel_{E^{l,p}}=\left\Arrowvert\left(I-\dfrac{d^2}{dv^2}\right)^{1/2}H\right\Arrowvert_{L_P}$.\\
		From the above assumptions, one can say that $H(\mathcal{L})$ is a weak type $(1,1)$ and bounded on some $L^p(M)$ spaces for all $p\in (1,\infty]$
	\end{coro}

	\subsection{Notations}
	Before the proof of the main results the following notations are very essential, also review of some supporting Lemma shall be considered
	\begin{defx}
		\ \\
		\normalfont
		A Lie group $G$ is said to be of polynomial growth if there exists a constant $K>0$, for some integer $D\geq1$ and $r\geq1$ such that
		$$K^{-1}r^D\leq V(r)\leq Kr^D$$
		This bounds automatically imply that $G$ is unimodular
	\end{defx}
	\noindent
	The third Lie group theorem states that, there exists a unique simply connected Lie group $G$ associated to Lie algebra $\mathfrak{g}$. Now suppose that $\left\{\tilde{V}_i\right\}_{i=1}^{\infty}$ form left invariant vector fields on the Lie group $G$ corresponding to the collection $\left\{V_i\right\}_{i=1}^{\infty}$. By the natural representation $\tilde{\sigma}$ of the Lie algebra $\mathfrak{g}$ it follows that $\tilde{\sigma}(\tilde{V}_i)=V_i$.\\
	By implication, we denote by $\sigma$ the representation of the Lie group $G$ in the space of bounded operators acting on $L^2(M)$ for which
	\begin{equation}\label{RRR}
	\sigma (g)f(v)=f(\sigma(g)v) \;\;\;\;\;\;\;\;\; \mbox{for all} \;\;\;\;\;\;\; f\in L^2(M).
	\end{equation}
	Observe that for any $f\in L^2(M)$ and $v\in M$
	$$\sigma(exp(s\tilde{V}_i))f(v)=f(exp(sV_i)v).$$
	Then, $d\sigma=\tilde{\sigma}$.\\
	Now for every function $\kappa\in L^1(M)$ the operator $\sigma(\kappa)$ can be define, by standard representation theory approach, as
	\begin{equation}\label{SSS}
	\sigma(\kappa)=\displaystyle\int_G\kappa(g)\sigma(g)dg
	\end{equation}
	where $g\in G$ and $\sigma(g)$ is an isometry acting on $L^p(M)$ spaces for $p\in[0,\infty]$ with respect to the positive Haar measure on Lie group $G$.\\
	The system of left invariant vector fields $\left\{\tilde{V}_i\right\}_{i=1}^{\infty}$ generates the Lie algebra $\mathfrak{g}$ and the control distance $\tilde{\rho}(g,q)$ related to the system cam be defined as
	$$\tilde{\rho}(g,q)=|g^{-1}q|, \;\;\;\;\; \mbox{where} \;\;\;\; |g|=\tilde{\rho}(g,u).$$
	Suppose that $\eta(s)$ form an admissible curve on the Lie group $G$ related to the vector fields $\left\{\tilde{V}_i\right\}_{i=1}^{\infty}$, for a function $f\in G$ one can write
	\begin{eqnarray*}
		\dfrac{d}{ds}f(\eta(s))&=& \beta_1(s)(\tilde{V}_1f)(\eta(s)) + \beta_2(s)(\tilde{V}_2f)(\eta(s)) + \cdots \\
		&=& \displaystyle\sum_{i=1}^{\infty}\beta_i(s)(\tilde{V}_if)(\eta(s)),
	\end{eqnarray*}
	with this, one can conclude that the graph $s\mapsto\sigma(\eta(s))v$, for every $v\in M$, form an admissible curve on $M$, for a function $w$ on $M$, we write
	$$\dfrac{d}{ds}w(\sigma(\eta(s))v)=\displaystyle\sum_{i=1}^{\infty}\beta_i(s)(V_iw)(\sigma(\eta(s)v)).$$
	The consequence of this is that
	\begin{equation}\label{TTT}
	\rho(v,\sigma(g)v)\leq|g|, \;\;\;\;\; \forall \; v\in M, \; g\in G.
	\end{equation}
	Like previous notation, $\tilde{B}(g,r)$ denotes the open ball related to $\tilde{\rho}$ with centre at $g$  and radius $r$.\\
	The volume $|\tilde{B}(g,r)|$ is independent of $g$ and one can have
	$$\tilde{U}(r)=|\tilde{B}(c,r)|=|\tilde{B}(g,r)|,$$
	where $c$ is the neutral element of $G$. Since all the Lie groups of polynomial growth satisfy $\eqref{KKK}$. Here, it implies that there exists a positive constant $K$ for which
	$$\tilde{U}(2r)\leq K\tilde{U}(r) \;\;\;\; \forall \;\; r>0.$$
	In what follows, we consider a function $\kappa\in L^1(G)$ and the related operator $\sigma(s)$.\\
	It can be ascertained that
	\begin{equation}\label{UUU}
	supp\; \sigma(\kappa)\subset \{(u,v)\in M^2:\rho(u,v)\leq R\}.
	\end{equation}
	For every open sets $V_i \subset M$, $m_i\in L^2(V_i,d\mu)$ and $R=\rho(V_1,V_2,\cdots)$ it follows that
	$$\langle\sigma(\kappa)m_1,m_2,\cdots\rangle=0.$$
	In the same manner we claim that $\sigma(\kappa)\geq0$ if
	\begin{equation}\label{VVV}
	\langle\sigma(\kappa)m_1,m_2,\cdots\rangle\geq 0
	\end{equation}
	for any $m_i\in L^2(M,d\mu)$ for which $m_i(v)\geq0$ for all $v\in M$.
	\begin{lem}\cite{3}
		\ \\
		\normalfont
		Let $\kappa\in L^1(M)$ and $supp\;\kappa \subset \tilde{B}(c,R)$. Then $\sigma(\kappa)$ satisfies the estimate $\eqref{UUU}$.\;\;\; $\Box$
	\end{lem}
	\noindent
		In what follows, the unit function on $M$ will be denoted by $\mathfrak{1}_M$, and is defined, for all $v\in M$, by $\mathfrak{1}_M(v)=1$
	
	\begin{lem}\cite{3}
		\ \\
		\normalfont
		Let $\kappa\in L^1(G)$ and let $\sigma(\kappa)$ be as defined in $\eqref{SSS}$ above. Hence
		\begin{equation*}
		\sigma(\kappa)\mathfrak{1}_M=\mathfrak{1}_M\displaystyle\int_G\kappa(g)dg.
		\end{equation*}
		In addition, for every function $\kappa\in L^1(G)$ if $\kappa(g)\geq0$ \;\; $\forall \;\; g\in G$, therefore
		$$\sigma(\kappa)\geq0$$
		which means that $\sigma(\kappa)$ admits $\eqref{VVV}$
	\end{lem}
	\begin{proof}
		\ \\
		\normalfont
		From $\eqref{RRR}$ and $\eqref{SSS}$, we have
		\begin{equation*}
		\sigma(\kappa)\mathfrak{1}_M(v)=\displaystyle\int_G\kappa(g)\mathfrak{1}_M(\sigma(g)v)dg=\displaystyle\int_G\kappa(g)dg.
		\end{equation*}
		Now, for any $v\in M$, if $m_i(v)\geq0$\; $(i=1,2,\cdots)$, we have
		$$\langle\sigma(g)m_1,m_2,\cdots\rangle\geq0.$$
		Therefore, for every $g\in G$, if $\kappa(g)\geq0$, we have
		$$\langle\sigma(\kappa)m_1,m_2,\cdots\rangle=\displaystyle\int_G\kappa(g)\langle\sigma(g)m_1,m_2,\cdots\rangle\geq0.$$
		Hence, the required result.
	\end{proof}
	\noindent
	Now, the operator $\tilde{\mathcal{L}}$ on the group $G$ will be defined as
	$$\tilde{\mathcal{L}}=-\displaystyle\sum_{i=1}^{\infty}\tilde{V}_i^2.$$
	The operator $\tilde{\mathcal{L}}$ spans semigroup acting on $L^p(G)$ space and the related heat kernel $\tilde{h}_s(g,q)=\tilde{h}_s(gq^{-1})$ satisfies $\eqref{MMM}$, which can be written as
	\begin{equation}\label{WWW}
	\dfrac{a_1}{\tilde{\mu}(B(g,\sqrt{s}))}e^{-\lambda^{\prime}|g|^2s^{-1}}\leq \tilde{h}_s(g)\leq\dfrac{a_2}{\tilde{\mu}(B(g,\sqrt{s}))}e^{-\lambda|g|^2s^{-1}}
	\end{equation}
	the estimate $\eqref{WWW}$ above is called Li-Yau inequality.
	The Poisson semigroup related to the operator $\tilde{\mathcal{L}}$ is
	$$\{e^{-s\tilde{\mathcal{L}}^{1/2}}\}_{s\geq0}.$$
	The convolution kernel related to the Poisson semigroup is denoted by $\tilde{q}_s(g)$ which can be written by the following integral, together with heat kernel, by subordinate formula:
	\begin{equation}\label{XXX}
	\tilde{q}_s(g)=\sigma^{-1/2}\displaystyle\int_0^{\infty} e^{-l}\tilde{h}_{s^2/(4l)}(g)\dfrac{dl}{l}.
	\end{equation}
	The next Lemma discussed the properties of the kernel $\tilde{q}_s$ which will be useful in the sequel.
	\begin{lem}\cite{3}\label{aaa}
		\ \\
		\normalfont
		Let $G$ be a simply connected Lie group of infinite dimensional, if $\tilde{q}_s(g)$ is the kernel related to Poisson semigroup and consider the operator $\tilde{\mathcal{L}}$ defined above. One have for all $s>0$
		\begin{equation}\label{YYY}
		\dfrac{a_1}{\tilde{U}(s+|g|)}\dfrac{s}{s+|g|}\leq \tilde{q}_s(g)\leq\dfrac{a_2}{\tilde{U}(s+|g|)}\dfrac{s}{s+|g|}
		\end{equation}
		By implication, there is a positive constant $K>0$ for which
		\begin{equation}\label{ZZZ}
		K^{-1}\tilde{q}_s(g)\leq \tilde{q}_{2s}(g)\leq K\tilde{q}_s(g)
		\end{equation}
	\end{lem}
	\noindent
	Let the characteristic function of the ball $\tilde{B}(c,s)$ be $\chi_{\tilde{B}(c,s)}$, then the next Lemma which is the complement of the Lemma $3.7$ follows.
	\begin{lem}\cite{3}
		\ \\
		\normalfont
		Keeping the conditions in the Lemma $3.7$, the following holds, for all $s>0,\;\;g\in G$
		\begin{equation}\label{AAA}
		\dfrac{\chi_{\tilde{B}(c,s)}}{\tilde{U}(s)}\leq K\tilde{q}_s(g)
		\end{equation}
		amd
		\begin{equation}\label{BBB}
		|\tilde{V}_i\tilde{q}_s(g)|\leq Ks^{-1}\tilde{q}_s(g)
		\end{equation}
	\end{lem}

	\subsection{Proof of the Main Results: Theorem 3.1}
	Fix
	$$h_s=\sigma(\tilde{h}_s)\;\;\;\;\;\; \mbox{and}\;\;\;\; q_s=\sigma(\tilde{q}_s)$$
	here, $h_s$ and $q_s$ represent the heat kernel of semigroup generated by $-\mathcal{L}$ and Poisson semigroup span by $-\mathcal{L}^{1/2}$ respectively, this is according to the standard representation theory. The operators $\exp(-s\mathcal{L}^{1/2})$ and $\exp(-s\mathcal{L})$ respectively have smooth kernels $q_s(u,v)$ and $h_s(u,v)$ for every $(u,v)\in M^2$ and $s>0$ (this is in line with H$\ddot{o}$rmander, 1967).\\
	And note
	$$V_ih_s=\sigma(\tilde{V}_i\tilde{h}_s)\;\;\;\;\;\; \mbox{and}\;\;\;\; V_iq_s=\sigma(\tilde{V}_i\tilde{q}_s)$$
	In the sequel, the following Lemmata will be useful.
	\begin{coro}[Dziuba$\acute{n}$ski and Sikora, 2019]
		\ \\
		\normalfont
		Suppose that $q_s$ is the kernel relating to the Poisson semigroup $\exp(-s\mathcal{L}^{1/2})$. Then there is a positive number $K>0$ for which
		\begin{equation}\label{CCC}
		K^{-1}q_s(u,v)\leq q_{2s}(u,v)\leq Kq_s(u,v)
		\end{equation}
		and
		\begin{equation}\label{DDD}
		|V_iq_s(u,v)|\leq Ks^{-1}q_s(u,v)
		\end{equation}
	\end{coro}
	
	\begin{lem}
		\ \\
		\normalfont
		Suppose that $q_s$ is as defined in Corollary 3.9. therefore there are positive constants $K$ and $a$ for which
		\begin{equation}
		q_s(u,v)\leq Kq_s(u^{\prime},v)\exp\left(a\dfrac{\rho(u,u^{\prime})}{s}\right)\;\;\; \mbox{for all}\;\; u,u^{\prime}, v\in\;M
		\end{equation}
	\end{lem}
	\begin{proof}
		\ \\
		\normalfont
		Define the function connecting $u$ and $u^{\prime}$, whose parameter is unit velocity, as $\eta:[0,\mathcal{T}]\longrightarrow M$. Now putting
		$$f(t)=q_s(\eta(t),v).$$
		Using $\eqref{DDD}$, we have
		$$|f'(t)|\leq Ks^{-1}f(t).$$
		By the reason of standard differential inequality,
		$$f(\mathcal{T})\leq f(0)\exp(a\mathcal{T}/s).$$
		For all the function $\eta$, taking minimum over $\mathcal{T}$ gives the result of the Lemma.
	\end{proof}
	
	For the proof of $\eqref{KKK}$, that is, the doubling condition.\\
	Fix, for all $r>0$
	\begin{equation}
	j_r(u,v)=\dfrac{\sigma(\chi_{\tilde{B}(c,r)})(u,v)}{\tilde{U}(r)}
	\end{equation}
	\noindent
	Applying equation $\eqref{AAA}$ together with Lemma $(3.6)$, we have
	\begin{equation*}
	0\leq j_r(u,v)\leq Kq_r(u,v)
	\end{equation*}
	Furthermore, with reference to Lemma $(3.5)$, we get
	\begin{equation*}
	\sup j_r(u,v)\subset B(v,r)
	\end{equation*}
	Consider a fix value $v\in M$ and $r>0$ with
	\begin{equation*}
	d_r=\underset{u\in B(v,r)}{\sup}j_r(u,v)
	\end{equation*}
	and
	\begin{equation*}
	D_r=\underset{u\in B(v,r)}{\sup}q_r(u,v)
	\end{equation*}
	As a matter of fact, $d_r\leq KD_r$\\
	Here it is enough to write
	\begin{equation}\label{Dd}
	d_r\geq U(v,r)^{-1} \;\;\;\; \forall\; v\in M\;\;\; \mbox{and}\;\; R>0
	\end{equation}
	By the condition of Lemma 3.6, it can be seen that
	\begin{equation*}
	1=\displaystyle\int_Mj_r(u,v)du\leq d_rU(u,r)
	\end{equation*}
	Using Lemma 3.10, yields
	\begin{equation*}
	q_r(u,v)\leq Ke^{2a}q_r(u^{\prime},v)
	\end{equation*}
	where $u,u^{\prime}$ are elements of $B(v,r)$\\
	Therefore
	\begin{equation*}
	e^{-2a}D_r\leq Kq_r(u^{\prime},v)
	\end{equation*}
	And so
	\begin{eqnarray*}
		e^{-2a}D_rU(v,r)&\leq& K\displaystyle\int_M q_r (u^{\prime},v)du^{\prime}\\
		&=& K\displaystyle\int_G \tilde{q}_r(g)dg=K.
	\end{eqnarray*}
	Hence
	\begin{equation}\label{Cc}
	e^{-2a}D_s\leq KU(v,s)^{-1}.
	\end{equation}
	Combining $\eqref{Dd}$ and $\eqref{Cc}$ gives
	\begin{equation}\label{EEE}
	U(v,s)^{-1}\leq d_s\leq D_s \leq KU(v,s)^{-1}
	\end{equation}
	Now, by estimate $\eqref{CCC}$ of Corollary 3.9, there is a constant $K>0$, for which
	\begin{equation*}
	K_{-1}q_s(u,v)\leq q_{2s}(u,v)\leq Kq_s(u,v).
	\end{equation*}
	Thus,
	\begin{equation}\label{FFF}
	D_s\sim D_{2s}.
	\end{equation}
	Using $\eqref{EEE}$ for $D_s$ and $D_{2s}$ together with $\eqref{FFF}$ results to doubling condition.\\
	Now, let proceed to the proof of $\eqref{MMM}$.\\
	Note that
	\begin{equation*}
	q_s(v,v)\leq D_s\leq KU(v,s)^{-1}.
	\end{equation*}
	Therefore
	\begin{equation*}
	\displaystyle\int|q_s(v,u)|^2=q_{2s}(v,v)\leq KU(v,s)^{-1}.
	\end{equation*}
	Now let
	\begin{equation*}
	D_{U_s}f(v)=U(v,s)f(v)
	\end{equation*}
	equivalently, the above can be written as
	\begin{equation*}
	\parallel D_{U_s} \exp(-sL^{1/2})\parallel_{1\rightarrow\infty}\leq K
	\end{equation*}
	for any positive value $s$. This means
	\begin{equation*}
	\parallel D_{U_s}\exp(-s^2L)\parallel_{1\rightarrow\infty}\leq \parallel D_{U_s}\exp(-s^2L^{1/2})\parallel_{1\rightarrow\infty}\times \parallel\exp(-s^2L)\exp(-s^2L^{1/2})\parallel_{1\rightarrow\infty}\leq K.
	\end{equation*}
	Hence,
	\begin{equation*}
	\parallel D_{U_s}V_i\exp(-s^2L)\parallel_{1\rightarrow\infty}\leq K\parallel D_{U_s}V_i\exp(-s^2L^{1/2})\parallel_{1\rightarrow\infty}\leq \dfrac{K}{s}.
	\end{equation*}

\end{document}